\documentclass[10pt]{article}
\newtheorem{theorem}{Theorem}

\title{Harmonic, Monogenic and Hypermonogenic Functions on Some Conformally Flat Manifolds in $R^n$ arising from Special Arithmetic Groups of the Vahlen Group}
\author{R.S. Krau{\ss}har \thanks{Department ofMathematical Analysis, Ghent University, Building S-22, Galglaan 2, B-9000 Ghent, Belgium. E-mail:{\tt krauss@cage.UGent.be}}\and John Ryan \thanks{Department of Mathematics, University of Arkansas, Fayetteville, AR 72701, USA {\tt jryan@uark.edu}}\and Qiao Yuying \thanks{Department of Mathematics, Hebai Normal University, Shijazhuang,  P. R. China}}
\date{~}
\begin{document}
\maketitle
\begin{abstract}\ This paper focuses on the development of harmonic and Clifford analysis techniques in the context of some conformally flat manifolds that arise from factoring out a simply-connected domain from $R^n$ by special arithmetic subgroups of the conformal group. Our discussion encompasses in particular the Hopf manifold  $S^1 \times S^{n-1}$, conformally flat cylinders and tori and some conformally flat manifolds of genus $g \ge 2$, such as $k$-handled tori and polycylinders. This paper provides a continuation as well as an extension of our previous two papers \cite{KraRyan1,KraRyan2}. In particular, we introduce a Cauchy integral formula for hypermonogenic functions on cylinders, tori and on half of the Hopf manifold. These are solutions to the Dirac-Hodge equation with respect to the hyperbolic metric. We further develop generalizations of the Mittag-Leffler theorem and the Laurent expansion theorem for cylindrical and toroidal  monogenic functions. The study of Hardy space decompositions on the Hopf manifold is also continued. Kerzman-Stein operators are introduced. Explicit formulas for the Szeg\"o kernel, the Bergman kernel and the Poisson kernel of half the Hopf manifold are given.  
\end{abstract}{\bf{Keywords}}: Conformally flat spin manifolds, automorphic forms, harmonic analysis, monogenic functions, hypermonogenic functions, boundary value problems, Hardy spaces. 

\section{Introduction}

\ Conformally flat manifolds are manifolds which have atlases whose transition functions are M\"obius transformations. From the viewpoint of Liouville's theorem (see, e.g., ~\cite{cn}) which states that the set of conformal maps in $R^n$ is reduced to that of M\"obius transformations for all $n \ge 3$, this class of manifolds can be regarded as higher dimensional generalizations of Riemann surfaces. Those that are further endowed with a spinor structure, are called conformally flat spin manifolds.  These are studied by many authors. See for instance ~\cite{GoKa,Chang} among others. 

\ One aspect of study in connection with this class of manifolds is the study of boundary value problems of PDE's and related operators.  For many classical two-dimensional boundary value problems, the complex Cauchy-Riemann operator pervades much of modern analytic techniques which further could be carried over to the context of two dimensional Riemann surfaces. 

\ For the treatment of their higher dimensional analogues, the less known Euclidean Dirac operator plays a similarly  powerful key role. As it is the linearization of the Euclidean Laplace operator, it can successfully be applied to a number of boundary value problems which arise from classical harmonic analysis, for instance the Dirichlet problem. We refer the interested reader for example to~\cite{lmq,lms,mcintosh,q,q1}. In recent decades, the term Clifford analysis has been used to describe the study of applications of Dirac operators in analysis. In particular, during the previous ten years much effort has been made to extend Clifford analysis from the context of Euclidean spaces to that of general Riemannian manifolds, see for example~\cite{calderbank,cn,mitrea}. Results in this direction include Cauchy type integral formulas and Plemelj projection formulas. However, in the general context of Riemannian manifolds little in the way of explicit formulas has been developed so far.  

\ For a number of examples of conformally flat spin manifolds, much  progress obtaining explicit formulas like Cauchy kernels and Green's formulas has been made very recently~\cite{v,lr,KraRyan1,KraRyan2,Ryan2003}. This paper provides a continuation as well as an extension of the recent works \cite{KraRyan1,KraRyan2,Ryan2003}. In particular, we introduce a Cauchy integral formula for hypermonogenic functions on cylinders, tori and on the Hopf manifold. Hypermonogenic functions are solutions to a Dirac-Hodge type equation with respect to the hyperbolic metric. 

\ We also develop generalizations of the Mittag-Leffler theorem and the Laurent expansion theorem for cylindrical and toroidal  monogenic functions. The study of Hardy space decompositions on the Hopf manifold is also extended. In particular, Kerzman-Stein operators are introduced and we give explicit formulas for the Szeg\"o kernel, the Bergman kernel and the Poisson kernel of half the Hopf manifold. We conclude by outlining some first rudiments for  a function theory on conformally flat manifolds of genus $g \ge 2$, such as Pretzel type manifolds and, more generally, $k$-handled tori and $k$-handled polycylinders.

\section{Preliminaries}

\  We start by introducing some basic notions and notations on real Clifford algebras over the Euclidean space $R^n$. For particular details we refer the interested reader for example to~\cite{p,bds}. Throughout this paper, $\{e_1,\ldots,e_n\}$ denotes the standard orthonormal basis in the Euclidean space $R^n$ and $Cl_n$ is its real Clifford algebra in which the relation $e_i e_j + e_j e_i = - 2 \delta_{ij}$ holds. This relation endows the Euclidian vector space $R^n$ with a multiplicative structure: Each vector $x \in R^n \backslash\{0\}$ has an inverse element, given by $-x/\|x\|^2$. The reversion anti-automorphism is defined by the relations $\widetilde{ab} = \tilde{b} \tilde{a}$, $\tilde{e_i} = e_i$ for all $i=1,2,\ldots,n$. The conjugation anit-automorphism is given by $\overline{ab} = \overline{b}\; \overline{a}, \overline{e_i} = -e_i$ for all $i=1,\ldots,n$. We further need the automorphism $*: Cl_n \rightarrow Cl_n$ which is  defined by the relations: $e_n^{*}=-e_n$, $e_i^{*}= e_i$ for $i=0,1,\ldots,n-1$ and $(ab)^{*}=a^{*}b^{*}$. 

\ Let $U \subseteq R^n$ be an open set. A real differentiable function $f: U \rightarrow Cl_n$ which satisfies $Df = 0$, resp. $fD = 0$, where $D := \frac{\partial }{\partial x_1} e_1 + \frac{\partial }{\partial x_2} e_2 + \cdots + \frac{\partial }{\partial x_n} e_n$  is the Euclidean Dirac operator, is called left (resp. right) monogenic or left (resp. right) Clifford-holomorphic. Due to non-commutativity, both classes of functions do not coincide with each other; however, for both classes an analogous function theory can be established, see e.g. \cite{bds} or elsewhere. The left and right fundamental solution to the $D$-operator is called the Euclidean Cauchy kernel which has the form $G(x-y) = \frac{x-y}{\|x-y\|^n}$. The Dirac operator factorizes the Euclidean Laplacian $\Delta = \sum_{j=1}^n \frac{\partial^2}{\partial x_j^2}$, viz $D^2 = - \Delta$ whose fundamental solution is the Green's kernel given by $H(x-y) = \frac{1}{\|x-y\|^{n-2}}$, whenever $n \ge 3$. Every real component of a monogenic function is hence harmonic. 

\ In contrast to complex analysis, the set of left (right) monogenic functions is not endowed with a ring structure. It forms only a right (left) $Cl_n$-module. In general, the composition of two monogenic functions does not result into a monogenic function, either. However, both the $D$-operator as well as the Laplacian is left quasi-invariant under M\"obius transformations. Following for example \cite{a,EGM87}, M\"obius transformations in $R^n$ can be represented in the way $T: R^n \cup \{\infty\} \rightarrow R^n \cup \{\infty\}$, $T(x)=(ax+b)(cx+d)^{-1}$ with coefficients $a,b,c,d$ from $Cl_n$ that can all be written as products of vectors from $R^n$ and that satisfy $a\tilde{d}-b\tilde{c} \in R^n \backslash\{0\}$ and $a^{-1}b, c^{-1}d \in R^n$ if $c \neq 0$ or $a \neq 0$, respectively. These conditions are often called Vahlen conditions. 

\ Now assume that $M(x)$ is a M\"obius transformation represented in the above mentioned form. If $f$ is a left monogenic (harmonic) function in the  variable $y=M(x)=(ax+b)(cx+d)^{-1}$, then, \cite{r85}, the function $J_1(M,x)f(M(x))$, resp. $J_2(M,x)f(M(x))$ is again left-monogenic (harmonic) in the variable $x$. Here, $J_1(M,x) = \frac{\widetilde{cx+d}}{\|cx+d\|^n}$ stands for the monogenic conformal weight factor and $J_2(M,x) = \frac{1}{\|cx+d\|^{n-2}}$ for the harmonic conformal weight factor. The set that consists of Clifford valued matrices $\left(\begin{array}{cc} a & b \\ c & d \end{array}\right)$ whose coefficients satisfy the above mentioned conditions, is called the general Vahlen group $GV(R^n)$ which generalizes $SL(2,C)$ to the higher dimensional context. The conformal invariance holds of course also for all subgroups of $GV(R^n)$. The particular subgroup $SV(R^n) =\{M \in GV(R^n)\;|\; a\tilde{d} - b \tilde{c} = 1\}$ is called the special Vahlen group and is purely generated by the inversion matrix and translation type matrices, as proved for instance in \cite{EGM87}. Its subgroup $SV(R^{n-1})$ has the special property that it acts transitively on the upper half-space  $H^{+}(R^n) = \{x \in R^n\;\ x_n > 0\}$. In view of this property, all M\"obius transformations associated to the group  $SV(R^{n-1})$ (as well as of all its subgroups) $T(x) = (ax+b)(cx+d)$ are well defined for each $x \in H^{+}(R^n)$ and $M(x)$ itself is an element from $H^{+}(R^n)$ for each $x \in H^{+}(R^n)$. Therefore, if $f$ is a function that is left monogenic on the whole half-space in the variable $x$, then $J_1(M,x)f(M(x))$ is so, too for all $x \in H^{+}(R^n)$. Similarly, for harmonic function, involving $J_2(M,x)$ instead. 

\ In this paper we also want to consider hypermonogenic functions. These are null-solutions to the hyperbolic Dirac operator $D_h: = x_n D f(x) + n Sc(-e_n f)(x) = 0$, where ($x_n \neq 0$).   The solutions to this system are called hyperbolic monogenic functions or often simply hypermonogenic functions for short.  One of the first to study this system extensively was H. Leutwiler in 1992~\cite{Leutwiler}. A Cauchy kernel function for hypermonogenic functions in the upper half-space has  been set up recently by  S.L. Eriksson in~\cite{E}.  Following, e.g., \cite{Hempfling}, also the hyperbolic D-operator shows a conformal invariance behavior: If $f$ is left hypermonogenic in the variable $y:=(ax+b)(cx+d)^{-1}$, where we assume that $\left(\begin{array}{cc} a & b \\ c & d \end{array}\right) \in SV(R^{n-1})$ then $F(x):=(cx+d)^{-1} f((ax+b)(cx+d)^{-1}) (x\tilde{c}+\tilde{d})$ turns out to be a left hypermonogenic in the variable $x$.  

\section{Monogenic, harmonic and hypermonogenic functions and some boundary value problems}

\ In this section we continue with the development of the rudiments of Clifford analysis, harmonic analysis and modified Clifford analysis on some  particular examples of conformally flat spin manifolds that arise from factoring out a simply-connected domain $U \subseteq R^n$ by discontinuous subgroups $\Gamma$ of the Vahlen group $GV(R^n)$. A general description of all discrete subgroups of the Vahlen group is given in \cite{EGM90} while in particular all congruence groups of the special hypercomplex modular group $\Gamma_p$ (its definition will be recalled in Section 3.3) are treated in \cite{Krapaper}. Readers who are not familiar with monogenic and harmonic spinor sections can find detailed information in \cite{r85}.  

\ Let $M$ be such a conformally flat manifold constructed by $M = U/\Gamma$, with $U,\Gamma$ as above. In all that follows $U\subseteq R^n$ is the universal covering space of the manifold $M$. As a consequence, there exists a well-defined projection map $p: U \rightarrow M, x \mapsto x\;(mod\; \Gamma)$. Let us write for each $x \in U$ for the image under the projection $x'=p(x)$. Consequently, for each open set $Q\subseteq R^n$ we write $Q'=p(Q)$. In this way $\Gamma$-invariant functions give rise to functions on $M$. Functions, that show an invariance behavior under the action of a discrete group are called automorphic functions. More generally one speaks of automorphic forms, when having a group invariance up to an elementary factor. 

\  For the classical theory of analytic automorphic forms in one complex variable, see for example \cite{schoeneberg}. For the basic theory of monogenic, harmonic and polymonogenic automorphic forms on general discrete subgroups of the Vahlen group we refer the interested reader to~\cite{Krahabil} and to~\cite{Krapaper}. Here in this paper, we explain how particular prototypes from this hypercomplex analytic type of automorphic forms can be used to tackle boundary value problems from Clifford and harmonic analysis on three particular families of spin manifolds that are constructed in this way. What follows provides an extension of our recent  papers~\cite{KraRyan1,KraRyan2} in which we produced some first fundamental steps in this direction. In particular we also develop some first results in the context of hypermonogenic functions.  

\subsection{Hopf Type Manifolds} 
{\bf Basic Background}.  Here, $U = R^n \backslash \{0\}$  and $\Gamma = \{m^k: k \in Z\}$ with an arbitrary but fixed $m \in N$ with $m > 1$. We consider one family of spinor bundles ${F_{1,m}}, m \in N$ over $S^1 \times S^{n-1}\cong U/\Gamma$, making the identification $(x,X) \Longleftrightarrow (m^k x,m^{k(n-1)/2}X)$ for every $k \in Z$ where $x \in R^n \backslash \{0\}$ and $X \in Cl_n$. By the identification $(x,X) \Longleftrightarrow (m^k x,m^{k(n-2)/2}X)$ a further family of spinor bundles, say $F_{2,m}$ is constructed. The first family gives rise to monogenic spinor sections and the second one with harmonic sections on these manifolds. Monogenic and harmonic functions on the associated Hopf type manifolds can be constructed from monogenic and harmonic automorphic forms that transform under the discrete dilatation group $\Gamma= \{m^k,l \in Z\}$ in the following way:

\begin{equation}
\label{dilatation}f(x) = m^{k(n-j)/2} f(m^k x) \quad\quad j=1,2.
\end{equation}

\ In~\cite{KraRyan2} explicit Cauchy (Green) formulae for such monogenic (harmonic) functions on the Hopf manifold have been developed. As these provide the starting point for the following investigation, we briefly quote them for convenience: 

\begin{theorem}\label{tcauchypr}
Let  $V'$ be a domain in $S^{1}\times S^{n-1}$. Suppose also that $S'$ is a strongly Lipschitz hypersurface in $V'$ and that $S'$ bounds a subdomain $W'$ and within $W'$ the domain $W'$ is contractable to a point. Then for each $y'\in W'$ and each left Clifford holomorphic section $f':V'\rightarrow F_{1,m}$ we have

\begin{equation}\label{cauchypr}f'(y')=\frac{1}{\omega_{n}}\int_{S'}{G'}_H(x',y')dp(n(x))f'(x')d\sigma'(x').
\end{equation}
Here, $\omega_n$ stands for the surface area of $S^n$, $n(x)$ the unit outer normal vector to $S:=p^{-1}(S')$ at $x:=p^{-1}(x')$ and $dp$ for the derivative of $p$. 
\end{theorem}

\ Here, ${G'}_H$ is the image of the function   $$G_H(x,y) = \sum\limits_{k=-\infty}^0 m^{k(n-1)/2} G(m^k x -y) + G(x) \Big[ \sum\limits_{k=1}^{\infty} m^{k(1-n)/2} G(m^{-k}x^{-1}-y^{-1})\Big] G(y),$$under the projection map $p$ which is exactly the $\Gamma$-periodization of the Euclidean Cauchy kernel function, transforming in the way (\ref{dilatation}).    

\ Under the same conditions, each harmonic section $f': V' \rightarrow F_{2,m}$ satisfies    

\begin{equation}
\label{greenpr}f'(y') = \frac{1}{\omega_n} \int\limits_{S'} \Big[ H_{H}(x',y') dp(n(x)) f'(x') + G_{H}(x',y') dp(n(x)) D_{H}f'(x')\Big] d\sigma'(x') \end{equation}

where $H_H$ is deduced from $G_H$ just by replacing the Euclidean Cauchy kernel function $G(x)$ by the harmonic one $H(x)$. 

\par\medskip\par {\bf Kerzman-Stein operators}.   

\ Let us next suppose that $S'$ is such a strongly Lipschitz hypersurface with these properties, dividing  the Hopf manifold $S^1 \times S^{n-1}$ into two disjoint domains, say $S'^{+}$ and $S'^{-}$. From the Cauchy (Green) formulas that we set up, one can now directly calculate Plemelj projection operators for $L^p(S') =\{\theta: S' \rightarrow F_j: \|\theta\|_{L_p} < \infty\}$ where $ 1 < p < \infty$. As shown in our previous work \cite{KraRyan2}, under these conditions one obtains the usual Hardy space decomposition of the form $L^p(S') = H^p(S'^{+}) \oplus H^p(S'^{-})$ in this context.  

\ This admits us now to carry over the standard Kerzman-Stein theory from the Euclidean context (see for instance \cite{BL} and elsewhere) to the framework of Hopf manifolds in an explicit way. To go into detail, let us denote the Cauchy-transform by ${\cal{H}}':=\frac{1}{2} I + C_{S'}$.  Adapting the calculations from \cite{BL} to the context of Hopf manifolds, one arrives at the  following explicit formula for the $L^p(S')$ adjoint of the Cauchy-transform$${\cal{H}}'\eta'(w') := \frac{1}{2} \eta'(w') + \frac{1}{\omega_n}P.V.\int\limits_{S'} \overline{dp(n(w))} {G'}_{H}(x',w')\eta'(x')d\sigma(x'). $$The adjoint Cauchy-transform ${\cal{H}}'\eta'(u')$ also turns out to be $L^p(S')$ bounded. Now it is natural to introduce the Kerzman-Stein operator on the Hopf manifold by $${\cal{A}}'\eta'(w'):=[{\cal{H}}'-{{\cal{H}}'}^{*}]\eta'(w') := \frac{1}{\omega_n} \int\limits_{S'} A(x',w') \eta'(x') d\sigma(x'),$$  introducing the Kerzman-Stein kernel as    $$A(x',w'):= {G'}_H(x',w')(dp(n(x))) - \overline{dp(n(w))} {G'}_H(x',w').$$  

\ Since $S'$ has a Lipschitz boundary, we thus can readily adapt the argument from \cite{BL}, Lemma 4.5, to the context of the Hopf manifold, and thus, we can establish that the operator ${\cal{I}}+{\cal{A}}'$ is invertible on $L^p(S')$. This has an important consequence. We can hence write the Szeg\"o projection ${\cal{P}}'$ on the space $L^p(S')$, $S'$ being a subdomain of the Hopf manifold with a strongly Lipschitz boundary, in the way $${\cal{P}}'={\cal{H}}'({\cal{I}}+{\cal{A}}'). $$  Notice that the Szeg\"o kernel is different from domain to domain, and it is very difficult to determine it explicitly for a given domain $S'$.  However, the Cauchy kernel is always the same, namely ${G'}_H$ on the Hopf manifold, independent of which domain $S'$ is considered. Since we have an explicit expression for ${G'}_H$ we can hence compute the Szeg\"o projection with this formula, without knowing the Szeg\"o kernel function. We can say more. In many cases, as for instance in cases when $S'$ is topologically {\it similar} to $p(\partial B)$, $B$ denoting a ball in $R^n$, this formula provides simpler numerical methods to determine the Szeg\"o kernel of $S'$, than to determine it directly, for example by applying a Gram-Schmidt algorithm which is numerically very instable.        

\par\medskip\par{\bf Szeg\"o kernel, Bergman kernel, and Poisson kernel of half the Hopf manifold}. 

\ Let us now focus on the particular decomposition, considering exactly half the Hopf manifold arising from the factorization $H^{+}(R^n)/\{m^k,k \in Z\}$ which shall be denoted by ${\cal{S}}_+$ in what follows. Hence we are in a particular situation of the general one considered so far.  For ${\cal{S}}_+$ we can readily establish a closed and explicit formula for the Szeg\"o kernel, and, also for the Poisson kernel and the Szeg\"o kernel; one can carry over the known formulas for the kernel functions of the upper half-space $H^{+}(R^n)$ \cite{cn} to the context considered here by using the explicit representation formulas for the Cauchy kernel of the Hopf manifold. The Szeg\"o kernel $S(x',w')$ of ${\cal{S}}_+$ thus turns out to be ${G'}_H(x',w') dp(x)|_{x=e_n}$ where $x' \in \partial {\cal{S}}_+$ and $w' \in {\cal{S}}_+$. The Poisson kernel is twice the real part of the Szeg\"o kernel: $P(x',w') = 2 Sc(S(x',w'))$. In fact this also follows readily from the Plemelj projection formulas that we developed in the previous subsection. This Poisson kernel then solves the Dirichlet problem for harmonic functions on ${\cal{S}}_+$ with given $L^p$ data ($1 < p < \infty$) on the boundary $\partial {\cal{S}}_+$. By similar arguments as presented in~\cite{cn}, we can further infer that the Bergman kernel $B(x',w')$ for $L^2$-integrable monogenic functions on ${\cal{S}}_+$ is given by $2 \frac{\partial }{\partial {x_n'}} S(x',w')$. Finally, the Bergman kernel for harmonic functions on ${\cal{S}}_+$ turns out to be the real part of the monogenic Bergman kernel function $B(x',w')$. 

\par\medskip\par{\bf Hypermonogenic kernels}. 

\ Let us now turn to also establish a Cauchy-integral formula for hypermonogenic functions on half the Hopf manifold ${\cal{S}}_+ = H^{R^n}/\{m^k,k\in Z\}$. To proceed in this direction we now need to work with the hypermonogenic kernels \cite{E}$$p_1(x,y):=\frac{(x-y)^{-1}}{\|x-y\|^{n-2} \|x-y^{*}\|^{n-2}}$$and$$p_2(x,y):=\frac{(x^{*}-y)^{-1}}{\|x - y\|^{n-2}\|x^{*}-y\|^{n-2}}$$ which provide the adequate analogous of the function $G(x-y)$ in the hyperbolic setting adapted to the hyperbolic metric $ds^2 = \frac{dx_1^2+\cdots+dx_n^2}{dx_n^2}$. In order to set up a Cauchy kernel for hypermonogenic functions on  ${\cal{S}}_+$ we next define the two basic hypermonogenic series $$h_1(x,y) = \sum\limits_{k=-\infty}^0 m^{k/2} p_1(m^k x ,y)  +  x^{-1} \Big[ \sum\limits_{k=1}^{\infty} m^{-k/2} p_1(m^{-k}x^{-1},y^{-1})\Big] y^{-1},$$and$$h_2(x,y) = \sum\limits_{k=-\infty}^0 m^{k/2} p_2(m^k x ,y)  +  x^{-1} \Big[ \sum\limits_{k=1}^{\infty} m^{-k/2} p_2(m^{-k}x^{-1},y^{-1})\Big] y^{-1}. $$

\ These then are the $\Gamma$-periodizations of the Euclidean hypermonogenic Cauchy kernel functions. They hence induce the hypermonogenic Cauchy kernel functions on the ${\cal{S}}_+$, simply by projecting them down.  This allows us to carry over the Cauchy integral formula from the standard setting established in \cite{E} to the context of the Hopf manifold in the following form 

\begin{theorem}
Suppose that $U' \subset {\cal{S}}_+$ is an open set and that $K' \subset U$ is a domain with a Lipschitz boundary $S'$. Let $f:U'\subset F_{h,m}$ be a left hypermonogenic function, where $F_{h,m}$ stands for the hypermonogenic section. Then   
\begin{eqnarray*}
f'(y') &=& \frac{2^{n-2}{y_n'}^{n-2}}{\omega_{n}} \int\limits_{S'} {h'}_1(x',y') dp(n(x)) f'(x') d\sigma(x) \\& - & \int\limits_{S'} {h'}_2(x',y') dp(n(x))^{*}  f^{*}(x') d\sigma(x')
\end{eqnarray*}
\end{theorem}

\ From this formula one can then depart to introduce analogous Plemelj projection operators and to establish a  Kerzman-Stein theory for hypermonogenic functions. This can be done in complete analogy to the monogenic case, so that we leave it to the reader as a simple exercise. 

\subsection{Conformally Flat Cylinders and Tori}
{\bf Basic Background}. Here, $U = R^n$ and $\Gamma$ is the discrete translation group  $$T(\Omega_k) = \Big\langle \left(\begin{array}{cc}  1 &   \omega_1\\ 0 & 1\end{array}\right), \cdots \left(\begin{array}{cc}  1 & \omega_k\\ 0 & 1\end{array}\right) \Big\rangle$$ where $\omega_1,\ldots,\omega_k$ are $R$-linear independent vectors $(k<n)$ spanning a $k$-dimensional lattice $Z \omega_1 + \cdots + Z \omega_k$. As explained in our previous work \cite{KraRyan2}, the decomposition of the lattice $\Omega_k$ into the direct sum of the sublattices $\Omega_l := Z {\omega_1} + \cdots + Z \omega_l$ and $\Omega_{k-l} := Z \omega_{l+1} + \cdots + Z \omega_{l}$ gives rise to  $k$ conformally inequivalent different spinor bundles $E^{(l)}$ on $C_k \cong R^n/T(\Omega_k)$ by making the identification $(x,X) \Longleftrightarrow (x+\underline{m}+\underline{n},(-1)^{m_1+\ldots+m_l} X)$ with $x \in R^n, X \in Cl_n$.  Monogenic, harmonic and hypermonogenic functions on the associated manifolds can be obtained by applying the projection map $p: R^n \rightarrow R^n/T(\Omega_k) \cong C_k$ on monogenic, resp. harmonic, resp. hypermonogenic functions that are automorphic functions on the associated discrete translation group $T(\Omega_k)$. Examples of monogenic and harmonic automorphic functions on $T(\Omega_k)$ are for instance the generalized Eisenstein type series whose basic properties are described in detail in \cite{Krahabil}, Chapter 2. They hence induce examples of monogenic and harmonic functions on cylinders and tori via the projection map.  In our previous works \cite{KraRyan1,KraRyan2} explicit formulas for the Cauchy kernel on $C_k$ in terms of particular variants of generalized Eisenstein type series have been established. Following ~\cite{KraRyan2}, for the cylinders $C_k$ with $k < n-1$  with the spinor bundles $E^{(l)}$, ($l \le k$), the Cauchy kernel functions are induced by the following variants of the generalized monogenic cotangent functions 

\begin{equation}
\cot_{1,k,l}(x,y) = \sum\limits_{\underline{m} \in \Omega_l,\underline{n} \in \Omega_{k-l}} (-1)^{m_1+\cdots+m_l} G(x+\underline{m}+\underline{n})\end{equation}and \begin{equation}\cot_{2,k,l}(x,y) = \sum\limits_{\underline{m} \in \Omega_l,\underline{n} \in \Omega_{k-l}} (-1)^{m_1+\cdots+m_l} H(x+\underline{m}+\underline{n}). 
\end{equation}

\ These variants of the monogenic cotangent series are the $T(\Omega_p)$-periodizations of the Euclidean Cauchy kernel, adapted according to the underlying spinor structure of the associated $C_k$.  Notice, these series only converge normally for $p<n-1$. In the case $k=n-1$ this construction can easily be adapted by coupling antipodal lattice points together. Due to the non-existence of $n$-fold periodic monogenic functions with only one point singularity of the order of the Cauchy kernel per period cell, one then needs to adapt this construction in the case $k=n$ in the way involving two point singularities $a$ and $b$ (incongruent modulo $\Omega_n$). For details, and the precise definition of the toroidal monogenic cotangent function $\cot_{1,k,l;a,b}(x,y)$ see \cite{KraRyan2}. If we replace in the representation of $\cot_{1,k,l}(x,y)$ the Euclidean kernel function $G$ by the harmonic kernel $H$,  then one gets the corresponding Green kernel functions, denoted by  $\cot_{2,k,l}(x,y)$ for cylindrical harmonic functions. However, notice that the convergence abscissa of the harmonic analogous of these monogenic cotangent type functions (replacing $G$ by $H$ in the series) is one less than in the monogenic case. This implies that the construction of the monogenic kernel function for the $n$-torus, only carries over to the $n-1$ cylinder in the harmonic context.  

\par\medskip\par  {\bf Harmonic Green kernels on $C_n$ expressed by Eisenstein type series}

\ A question that remained unsolved in our previous two papers is to ask whether these constructions can further be extended to the case  $k=n$ in the context of harmonic functions. In this paper we finally set up a solution to this problem. To get a harmonic $n$-fold periodic cotangent type series involving only poles of the order of the Euclidean Cauchy kernel, one first might think of trying a further coupling of antipodal lattice points. This however does not work here so easily, since the Euclidean harmonic Green's kernel is an even function and not an odd one. Further minus signs hence cannot be generated by a further coupling of antipodal points, since they are eliminated by the modulus. The solution to this problem is to take for mutually incongruent modulo $\Omega_n$ points $a_1,a_2,a_3,a_4$ from $R^n \backslash \Omega_n$ and to make the following construction 

\begin{eqnarray*}& & \cot_{1,n,l;a_1,a_2,a_3,a_4}(x,y) \\& := & H(x-a_1) + H(x-a_2) + H(x-a_3) + H(x-a_4)\\ & & +   \sum\limits_{\underline{m} \in \Lambda_l, \underline{n} \in \Lambda_{k-l}} (-1)^{m_1+\cdots+m_l} \Big[       H(x-a_1+\underline{m}+\underline{n}) + H(-a_1-\underline{m}-\underline{n}) \\& &  +H(x-a_2+\underline{m}+\underline{n}) + H(-a_2-\underline{m}-\underline{n}) \\& &  -H(x+a_3+\underline{m}+ \underline{n})- H(-a_3 - \underline{m} - \underline{n}) \\& &  -H(x+a_4+\underline{m}- \underline{n})- H(-a_4 - \underline{m} - \underline{n})\Big]. 
\end{eqnarray*}

\ Here, $\Omega_l^{+} = N \omega_1 + Z \omega_2 + \cdots + Z \omega_l$ $\cup$ $N \omega_2 + Z \omega_3 + \cdots + Z \omega_l$ $ \cup \cdots \cup N \omega_l$ is the positive part of the lattice $\Omega_l$, according to the direct lattice decomposition as described in~\cite{Krahabil}, Chapter 2. Its projection $p(\cot_{1,n,l;a_1,a_2,a_3,a_4}(x,y)) := {\cot'}_{1,n,l;a_1,a_2,a_3,a_4}(x',y')$ induces then a local Green kernel on the torus which has four isolated point singularities and one can set up a local Green's formula for the torus. \begin{theorem}Suppose that $U'$ is a domain in $C_{n}$ and that $f':U'\rightarrow Cl_{n}$ is annihilated by the toriodal Laplace  operator $D'=p(\Delta)$ in $U'$. Suppose also that $V'$ is a subdomain of $U'$ whose closure also lies in $U'$ and $V'$ has a strongly Lipschitz boundary. Furthermore $a_2,a_3,a_4$ do not belong to the closure of $V'_{a_1}$. Then for each $y'\in V'$

\begin{equation}
\label{inttorus}f'(y')=\frac{1}{\omega_{n}}\int_{\partial V'_{a_1}}\cot'_{1,n,l;a_1,a_2,a_3,a_4}(x',y'_{a_1})(dp_{n}n(x))F'(x')d\sigma'(x')
\end{equation}
where $y'_{a_1}=p(y-a_1)$.
\end{theorem} 
{\bf Mittag-Leffler and Laurent expansion theorems for cylindrical and toroidal functions}. 

\ Next we give a generalization of Mittag-Leffler's theorem and the Laurent expansion theorem in the context of cylinders and tori. To this end, we now have additionally to also consider the partial derivatives of the cotangent functions. According to~\cite{Krahabil} Chapter 2.1, these partial derivatives represent the higher dimensional monogenic and harmonic analogous of the translative Eisenstein series of higher pole order. We shall see that they will play the analogous role of the partial derivatives of the Euclidean Cauchy kernel in Mittag-Leffler type and Laurent expansion type theorems in the context of conformally flat cylinders and tori. Let us consider slightly more general the following variants  $$\epsilon_{1,k,l,{\bf m}}(x):= \frac{\partial^{m_1+\cdots+m_n}}{\partial x_1^{m_1} \cdots \partial x_n^{m_n} }\cot_{1,k,l}(x,0)$$and$$\epsilon_{2,k,l,{\bf m}}(x):= \frac{\partial^{m_1+\cdots+m_n}}{\partial x_1^{m_1} \cdots \partial x_n^{m_n} }\cot_{2,k,l}(x,0),$$adapted to the underlying spinor structure. These functions allow us to adapt the standard Mittag-Leffler theorem from the Euclidean setting to the context of cylinders and tori with the particular spinor bundles under consideration in the following form:

\begin{theorem}Let $k \in \{1,\ldots,k\}$ Let ${a'}_1,{a'}_2,\ldots,{a'}_p \in C_k$ be a finite set of points. Suppose that $f': C_k \backslash\{{a'}_1,\ldots,{a'}_k\} \rightarrow E_j$ be a cylindrical/toroidal left monogenic (harmonic) function which has atmost isolated poles at the points ${a'}_i$ of the order $K_i$. Then there exists a cylindrical/toroidal left monogenic (harmonic) function $\phi: C_k \rightarrow E_j$ such that 
\begin{equation}
f'(x')= \sum\limits_{i=1}^p \sum\limits_{m=0}^{K_i-(n-1)}\sum\limits_{m=m_1+m_2+\cdots+m_n} \Bigg[{\epsilon'}_{1,k,l,{\bf m}}(x'-{a'}_i)\Bigg] + \phi'(x')
\end{equation}
in the monogenic case, or 
\begin{equation}
f'(x')= \sum\limits_{i=1}^p \sum\limits_{m=0}^{K_i-(n-1)}\sum\limits_{m=m_1+m_2+\cdots+m_n} \Bigg[{\epsilon'}_{2,k,l,{\bf m}}(x'-{a'}_i)\Bigg] + \phi'(x')
\end{equation}
in the harmonic case, respectively. In the case $k=n$ and ${\bf m} = 0$ one takes $\cot_{1,n,l;a_i,b}(x,0)$ where $b$ can be chosen arbitrarily such that $b \not\in \Omega_n$ and $b \not\equiv a$ (mod) $\Omega_n$. Similary, one proceeds in the harmonic case for $k=n-1$ with ${\bf m}=0$, and in the case $k=n$ with ${\bf m}=0$ one takes $ \cot_{1,n,l;a_i,b_2,b_3,b_4}(x,0)$, chosing $b_2,b_3,b_4$ adequately.  
\end{theorem}

\ Remark: In the particular case of the $n$-torus $C_n$ a special effect appears in this context, which does not appear in the Euclidean case, neither in case of the other $k$-cylinders with $k < n$.  In the particular case of the torus, the function $\phi'(x')$ always reduces to a single constant. This is a consequence of the fact that there are no non-constant entire monogenic (harmonic) $n$-fold periodic functions in $R^n$. A left toroidal monogenic (harmonic) function having at most isolated singularities on the torus can thus be represented up to a Clifford constant by a finite sum of the cotangent type series and a finite amount of their partial derivatives. 

\ Similarly, one can also establish a Laurent expansion theorem for cylindrical/toroidal left monogenic (harmonic) functions. Their local behavior near an isolated singularity $a' \in C_k$ can be described in the form $$f'(x') = \sum\limits_{|{\bf m}| \ge 0} {\epsilon'}_{j,k,l,{\bf m}}(x'-{a'}) +\phi'(x'-a'),$$where $\phi'(x'-a')$ is a cylindrical/toroidal left monogenic (harmonic) function in a sufficiently small neighborhood around $x'=a'$, satisfying $\lim_{x'\rightarrow a'} \phi'(x'-a') = \phi'(0)$. In the particular case when there are no further singularities on $C_k$, this representation holds globally on $C_k$, whence we are in a special case of the previous theorem. However, we have to keep in mind, that a left monogenic toroidal function must at least have two poles of order $n-1$ or at least a pole of higher order. 

\par\medskip\par{\bf Hypermonogenic kernels}. 

\ Next we want to establish integral formulas for hypermonogenic functions on the cylinders $C_k$ with $k < n$. To proceed in this direction we need to introduce hypermonogenic generalizations of the cotangent functions $\cot_{1,k,l}$. In view of the conformal invariance of $D_{hyp}$, a hypermonogenic function $f$ remains monogenic under all  translations of the form $f(z+\omega)$ whenever $\omega \in R^{n-1}$. If $\Omega_{k}$ is thus a lattice with the property $Sc(-e_n \Omega_k) = 0$, then all transformations of the form $f(z+\omega)$ with such lattice elements remain hypermonogenic. Due to the condition $Sc(-e_n \Omega_k) = 0$ we need to restrict in this context to consider only $p$-dimensional lattices with $p < n$. In view of the conformal invariance and Weierstra{\ss} convergence theorem, it is thus natural to introduce the hypermonogenic versions of the $\cot_{1,k,l}$ functions as follows, $$c_1(x,y)= \sum\limits_{\underline{m} \in \Omega_l^+,\underline{n} \in \Omega_{k-l}^-} (-1)^{m_1+\cdots+m_l} p_1(x+\underline{m}+\underline{n},y)$$and$$c_2(x,y) = \sum\limits_{\underline{m} \in \Omega_l^+,\underline{n} \in \Omega_{k-l}^-} (-1)^{m_1+\cdots+m_l} p_2(x+\underline{m}+\underline{n},y),$$where we assume $\Omega_p$ to be a $p$-dimensional lattice that is completely contained in $span_{R}\{e_1,\ldots,e_{n-1}\}$. Due to the additional factor $\|x-y^{*}\|^{n-2}$ that appears in the denominator of $p_1(x,y)$ and $p_2(x,y)$, we get better convergence conditions than in the monogenic case. Therefore, it is not necessary to add convergence preserving terms in the case $p=n-1$ which was necessary in the monogenic case. For the $k$-cylinders with $k < n$ we thus arrive at the following Cauchy integral formula: 

\begin{theorem}
Let $\Omega_k$ be a $k$-lattice that lies completely in $span_R\{e_1,\ldots,e_{n-1}\}$. Let $C_k$ be the cylinder induced by $R^n/\Omega_k$. Suppose that $U' \subset C_k$ is an open set and that $K' \subset U$ is a domain with a Lipschitz boundary $S'$. Let $E_j$ be one of the above mentioned spinor sections and suppose that $f:U'\subset E_j$ is a left hypermonogenic function. Then   
\begin{eqnarray*}
f'(y') &=& \frac{2^{n-2}{y_n'}^{n-2}}{\omega_{n}} \int\limits_{S'} {c'}_1(x',y') dp(n(x)) f'(x') d\sigma(x) \\& - & \int\limits_{S'} {c'}_2(x',y') dp(n(x))^{*}  f^{*}(x') d\sigma(x')
\end{eqnarray*}
\end{theorem}

\ Remark: Notice that the operator $D_{h}$ is not conformally invariant under translations in the $e_n$ direction. This construction hence cannot be carried over directly to the context of the $n$-torus $C_n$. This is a significant difference to the monogenic or Euclidean harmonic case. 

\subsection{Manifolds of higher genus}

\ All these examples are manifolds with a genus being at most one.  We conclude this paper by also dedicating a few words on conformally flat spin manifolds of genus $g \ge 2$.  For the construction of examples of conformally flat spin manifolds of higher genus, such as pretzel-type manifolds, or, more generally $k$-handled tori  and $k$-handled polycylinders, we can use special arithmetic congruence groups of finite index in the special hypercomplex modular group $$\Gamma_p := \langle T_{e_1},\ldots,T_{e_p}, J\}.$$such as for example its homogeneous principal congruence subgroups of level $N$ $$ \Gamma_p[N] = \Bigg\{ \left( \begin{array}{cc} a & b \\ c & d\end{array} \right) \in \Gamma_p\;\Bigg|\; a-1,b,c,d-1 \in N {\cal{O}}_p\Bigg\} $$or the following basic examples of congruence groups with conductor $N$ $$ {\Gamma_p}_0[N] = \Bigg\{ \left( \begin{array}{cc} a & b \\ c & d\end{array} \right) \in \Gamma_p\;\Bigg|\;  c \in N {\cal{O}}_p\Bigg\} $$which satisfy $\Gamma_p[N] \subseteq {\Gamma_p}_0[N] \subseteq \Gamma_p.$ Here, ${\cal{O}}_k = \sum_{A \in P(1,\ldots,k)} Z e_A$ denotes the standard order in the real Clifford algebra $Cl_n$.  Factoring out the upper half-space $H^{+}(R^n)$ by groups of the type $\Gamma_p[N]$ or ${\Gamma_p}_{0}[N]$ leads to examples that belong to this family of conformally flat manifolds with genus $g \ge 2$ when choosing $N$ sufficiently large. Further examples of conformally flat  manifolds of higher genus can be obtained by factoring Cartesian products of half-spaces by groups of the type $\Gamma_p[N]$ or ${\Gamma_p}_{0}[N]$. This construction method provides a counterpart to the gluing method discussed in the recent paper~\cite{Ryan2003}.

\par\medskip\par  Monogenic and harmonic automorphic forms on these groups hence give rise to monogenic and harmonic automorphic forms on these families of manifolds. Some non-trivial examples of (poly-)monogenic automorphic forms for congruence groups of the hypercomplex modular group $\Gamma_p$ which include these groups, have recently been constructed in \cite{Krahabil,Krapaper} on some symmetric spaces which provide generalizations of the analytic Eisenstein- and Poincar\'e series.  Let us now analyze which of them can be used (and how) in order to define non-trivial monogenic and harmonic functions on $H^{+}/\Gamma$ where $\Gamma$ is either a group belonging to the families $\Gamma_p[N]$ or ${\Gamma_p}_0[N]$. 
\par\medskip\par

\ Following \cite{Krahabil}, the function series 

\begin{equation}
\label{eisen1}{\cal{G}}^{p,N}(x) = \sum\limits_{M:{\cal{T}}_p[N] \backslash \Gamma_p[N]} J_1(M,x),
\end{equation}
where $M:T_p[N] \backslash \Gamma_p[N]$ means that the matrices $M$ run through a system of representatives of the right cosets of $\Gamma_p[N]$ modulo its subgroup of translation matrices denoted by ${\cal{T}}_p[N]$, converges normally on $H^{+}(R^n)$ under the condition $p < n-1$; it satisfies$${\cal{G}}^{p,N}(x) = J_1(M,x){\cal{G}}^{p,N}((ax+b)(cx+d)^{-1}) \quad\quad \forall \left(\begin{array}{cc} a & b \\ c & d\end{array}\right) \in \Gamma_p[N] $$and its limit towards $+e_n \infty$ equals $1$ under the condition $N \ge 3$. Hence it is a non-trivial monogenic automorphic function for $\Gamma_p[N]$ on $H^{+}$ under these conditions, and induces, again under these conditions a left monogenic (non-vanishing) function on the associated spin manifold $H^{+}/\Gamma_p[N]$ which is of genus $g \ge 2$, provided $N$ is chosen sufficiently large.   However, as explained in \cite{Krahabil,Krapaper}, for all larger groups, including $\Gamma_p[1], \Gamma_p[2]$ and also for the other basic congruence groups ${\Gamma_p}_0[N]$ the construction in (\ref{eisen1}) vanishes identically. This is basically due to the fact that the negative identity matrix is then included. Such series do hence not provide us with monogenic functions on manifolds that are parameterized by these larger groups. This reflects a fundamental algebraic property. The groups $\Gamma_p[N]$ with $N \ge 3$ are the largest congruence groups in $\Gamma_p$ with the property that they do not contain the negative identity matrix. All of them are normal subgroups in $\Gamma_p$ which is not the case for all other congruence groups lying in between $\Gamma_p[N]$ and $\Gamma_p$. 

\par\medskip\par

\ In the harmonic case, the situation is slightly different. Due to the fact that $H(x)$ is an even function which takes only non-negative positive values, the series  
\begin{equation}
\label{eisen2}{\cal{H}}^{p,N}(x) = \sum\limits_{M:{\cal{T}}_p[N] \backslash \Gamma_p[N]} J_2(M,x)
\end{equation}
do also not vanish for $N=1,2$. As explained in \cite{Krahabil}, their limit towards $+e_n \infty$ equals $+2$ or $2^{p+1}$ in case $N=1$ or $N=2$, respectively. Hence the projection $p: H^{+}(R^n) \rightarrow H^{+}(R^n)/\Gamma_p[N]: p(x) ={\cal{H'}}^{p,N}(x')$ induces a well defined non-constant harmonic function on $H^{+}(R^n)/\Gamma_p[N]$ for each $N \ge 1$ where $p < n-2$. However, one should mention that the two additional cases $N=1,2$ do not provide a significant enrichment within this context of spin manifolds, since both $H^{+}(R^n)/\Gamma_p[1]$ and $H^{+}(R^n)/\Gamma_p[2]$ have genus smaller than $2$. However, one gets more in the harmonic case in the following sense. The crucial point is that in the harmonic case also all the series

\begin{equation}
\label{eisen3}{\cal{H_0}}^{p,N}(x) = \sum\limits_{M:{{\cal{T}}_p}_0[N] \backslash {\Gamma_p}_0[N]} J_2(M,x)
\end{equation}  
are non-vanishing automorphic forms on the basic congruence groups ${\Gamma_p}_0[N]$ for any $N \in N$ which is not the case in the monogenic case, since the negative identity matrix is included in all these groups. 

\ However, due to the evenness of the harmonic weight function, one also gets a positive limit towards $+e_n \infty$ for all these groups. The projection $H^{+}(R^n) \rightarrow H^{+}(R^n)/{\Gamma_p}_0[N]: p(x) ={\cal{H_0'}}^{p,N}(x')$ induce hence for all $N \ge 1$ non-trivial harmonic functions on $ H^{+}(R^n)/{\Gamma_p}_0[N]$ which provide us with further examples of conformally flat manifolds of genus $g \ge 2$, whenever $N$ is chosen sufficiently large. 

\ Two natural questions arise immediately in this context. First, is it possible to also construct monogenic functions on the manifolds $H^{+}/{\Gamma_p}_0[N]$? Second, is it possible to construct monogenic and harmonic functions on those manifolds of this type attached to the parameter $p=n-1$ and for $p=n-2,n-1$, respectively? Notice that the convergence abscissa of the series (\ref{eisen1}) is only $p < n-1$ in the monogenic case and $p< n-2$ in the harmonic case, for all groups under consideration. In \cite{Krahabil,Krapaper} a construction has been proposed to also get monogenic and harmonic automorphic forms on $\Gamma_{n-1}$ and all its congruence subgroups $\Lambda$  that satisfy $\Gamma_p[N] \subseteq \Lambda \subseteq\Gamma_p$ for an $N \in {\bf N}$:       
\begin{equation}
{\cal{E}}^{\Lambda}_1(x,y) = \sum\limits_{M: T(\Lambda) \backslash \Lambda} \overline{\widetilde{J_1(M,x)}} \widetilde{J_1(M,y)}
\end{equation}

and

\begin{equation}{\cal{E}}^{\Lambda}_2(x,y) = \sum\limits_{M: T(\Lambda) \backslash \Lambda} \overline{\widetilde{J_2(M,x)}} \widetilde{J_2(M,y)}
\end{equation}

\ Here, $T(\Lambda)$ stands for the subgroup of translation matrices that is contained in $\Lambda$. This type of series converges in the monogenic case even for $p=n-1$ whenever $n \ge 3$. Similarly, in the harmonic case for $n \ge 4$, due to having now two automorphy factors. Because of the two automorphy factors one furthermore gets a non identically vanishing series for all $\Gamma_p$ and all their congruence subgroups, including in particular all the groups $\Gamma_p[N]$ and ${\Gamma_p}_0[N]$. The limit towards $+e_n \infty$ turned out to be always positive. Indeed, this series is an automorphic form on these groups; however, its projection does not induce a function on the families of manifolds of the type $H^{+}(R^n)/\Gamma_p[N]$ nor  $H^{+}(R^n)/{\Gamma_p}_0[N]$, respectively, since it is a function in two vector variables defined on the Cartesian product of two half-spaces! It hence induces examples of monogenic and harmonic functions on another type of conformally flat spin manifolds, namely on those that arise from factoring out $H^{+}(R^n) \oplus H^{+}(R^n)$ by ${\Gamma_p}[N]$ and ${\Gamma_p}_0[N]$. This is a different type of manifold, however, choosing $N$ sufficiently large, these will also give examples of conformally flat manifolds of higher genus, although not topologically equivalent to the others treated before.

\  Many interesting questions are thus still open in this context. For a future development and a more quantative understanding it would be extremely important to find an explicit relation between the conductor of the particular congruence groups $\Gamma_p[N]$ and ${\Gamma_p}_0[N]$ and the genus of the manifold. In the classical two dimensional case, such an explicit relation is given in terms of the genus formula \cite{schoeneberg} which can be deduced as a special case of Euler's polyhedra formula. This problem is closely connected to the Riemann Roch theorem. In the context considered here, we would need a higher dimensional analogue of this machinery. A further task  of central importance is to find in particular those monogenic and harmonic functions on these manifolds that serve as Cauchy kernels or Green kernels.  Maybe the above mentioned series could provide us some useful building blocks for the construction of such a kernel, since they have the good property of being the $\Gamma_p[N]$-periodizations (${\Gamma_p}_0[N]$-periodizations) of the fundamental solution of the Euclidean Dirac operator.

\end{document}